\documentclass[11pt]{article}
\usepackage{iftex}
\ifpdf\else
  \ifXeTeX\else
    \PassOptionsToPackage{dvipdfmx}{graphicx}
    \PassOptionsToPackage{dvipdfmx}{xcolor}
    \PassOptionsToPackage{dvipdfmx}{hyperref}
    
  \fi
\fi
\usepackage[T1]{fontenc}
\usepackage[utf8]{inputenc}
\usepackage{lmodern}
\usepackage{amsmath,amssymb,amsthm,mathtools}
\usepackage{aliascnt}
\usepackage{authblk}
\usepackage{tikz,graphicx}
\usetikzlibrary{calc}
\usepackage[margin=1in]{geometry}
\usepackage{booktabs,array,microtype}
\usepackage{enumitem}
\usepackage[hidelinks]{hyperref}
\usepackage[capitalise,nameinlink,noabbrev]{cleveref}

\newtheorem{theorem}{Theorem} 
\newtheorem{lemma}[theorem]{Lemma}
\newtheorem{conjecture}[theorem]{Conjecture}

\newtheorem{claim}{Claim}
\crefname{claim}{Claim}{Claims}
\Crefname{claim}{Claim}{Claims}

\theoremstyle{definition}
\newtheorem{definition}{Definition}

\newcommand{\cP}{\mathcal P}

\AtBeginEnvironment{proof}{\setcounter{case}{0}}
{\setlength{\leftmargini}{1.5\parindent}
	\begin{itemize}
		\setlength{\itemsep}{-1.1mm}}
	{\end{itemize}}

\begin{document}
\title{A sharp minimum degree condition for Hamiltonicity in line graphs of 3-hypergraphs
\thanks{Kenta Ozeki was supported by JSPS KAKENHI, Grant Numbers 22K19773, 23K03195, 26K00616 and 26K00618. Leilei Zhang's work was supported by  the National Natural Science Foundation of China (Grant Nos. 12601664, 12571365), JSPS KAKENHI (Grant No.~25KF0036), the Hubei Provincial Natural Science Foundation of China (Grant No.~2025AFB309), the China Postdoctoral Science Foundation (Grant No.~2025M773113), the Fundamental Research Funds for the Central Universities, Central China Normal University (Grant No.~CCNU24XJ026).}}
	
\author[1]{Kenta Ozeki \footnote{Email: \texttt{ozeki-kenta-xr@ynu.ac.jp.}}}
\author[1,2]{Leilei Zhang\footnote{Email: \texttt{leilei-zhang-gt@outlook.com.}}}

\affil[1]{\small Faculty of Environment and Information Sciences, Yokohama National University, Yokohama, Japan }
\affil[2]{\small School of Mathematics and Statistics, Key Laboratory of Nonlinear Analysis \& Applications (Ministry of Education), Central China Normal University, Wuhan, China}
\date{}
\maketitle
\vspace{-2.5em}
\begin{abstract}
 
Let $G$ be the line graph of a $3$-hypergraph, and let $n$ be its order. We prove that if $G$ is $3$-connected, $n\ge4936$, and $\delta(G)\ge(n+5)/10$,  then $G$ is Hamiltonian if and only if it does not belong to an explicitly described family of exceptional graphs.
This extends the sharp minimum degree bound of Lai, Shao, and Zhan for $3$-connected claw-free graphs beyond the claw-free setting, since line graphs of $3$-hypergraphs may contain induced claws. The coefficient $1/10$ cannot be decreased, and the constant
$5$ in the numerator cannot be replaced by $4$.
\end{abstract}
	
\section{Introduction}
	
All graphs are finite and simple unless explicitly called multigraphs. Multigraphs may have parallel edges but no loops. We use the standard notation and terminology of Bondy and Murty~\cite{BM}.

The study of minimum degree conditions for Hamiltonicity is a classical topic in graph theory. 
Dirac~\cite{D} proved that every graph of order $n\ge 3$ with minimum degree at least $n/2$ is Hamiltonian. A graph is called \emph{claw-free} if it contains no induced subgraph isomorphic to $K_{1,3}$. For claw-free graphs, weaker minimum degree conditions suffice under suitable connectivity assumptions. Matthews and Sumner~\cite{MS} proved that every $2$-connected claw-free graph of order $n$ with minimum degree at least $(n-2)/3$ is Hamiltonian. For regular graphs, Li~\cite{ML1} proved that every $2$-connected $k$-regular claw-free graph of order $n$ is Hamiltonian, provided $k \ge (n+2)/4$.

For $3$-connected claw-free graphs, several stronger results have been obtained. Let $G$  be a $3$-connected claw-free graph of order $n$ and minimum degree $\delta$. Zhang~\cite{Z} proved that $\delta\ge (n-3)/4$ guarantees Hamiltonicity, while Li~\cite{ML} established the sufficient condition $\delta\ge (n+5)/5$. Li, Lu, and Liu~\cite{LLL} further proved that 
$\delta \ge (n+7)/6$ suffices.

To state the results below, we recall the closure operation introduced by Ryj\'a\v{c}ek~\cite{R} in 1997.
The \emph{closure} $\operatorname{cl}(G)$ of a claw-free graph $G$ is obtained by repeatedly completing the neighborhood of an eligible vertex, where a vertex is \emph{eligible} if its neighborhood is connected but not complete, until no eligible vertex remains. Ryj\'{a}\v{c}ek \cite{R} proved that $\operatorname{cl}(G)$ is uniquely determined, is the line graph of a triangle-free graph, and is Hamiltonian if and only if $G$ is Hamiltonian.

Using the closure operation, Kov\'a\v{r}\'ik,  Mula\v{c} and Ryj\'a\v{c}ek~\cite{KMR} proved that every $2$-connected claw-free simple graph $G$ of order $n\ge 153$ with $\delta(G)\ge (n+39)/8$ is either Hamiltonian or has its closure in one of five explicitly described families of graphs. 
For $3$-connected claw-free graphs of order $n$ and minimum
degree $\delta$, Favaron and Fraisse~\cite{FF} established
the sufficient condition $\delta\ge(n+38)/10$ for all $n$,
where the coefficient $1/10$ is best possible.
Lai, Shao, and Zhan~\cite{LSZ} lowered the bound further and characterized the exceptional graphs in the following theorem.
	
\begin{theorem}[Lai, Shao, and Zhan~\cite{LSZ}]
\label{LSZ_thm}
Let $G$ be a $3$-connected claw-free simple graph of order $n\ge 196$ with $\delta(G)\ge (n+5)/10$. Then either $G$ is Hamiltonian or $\operatorname{cl}(G)=L(H)$, where $H$ is obtained from the Petersen graph $P$ by attaching $(n-15)/10$ pendant edges to each vertex of $P$.
\end{theorem}

For further results in this direction, see~\cite{CLX, ML2}. 

Although Theorem \ref{LSZ_thm} is stated for claw-free graphs, its essential content lies in the case of ordinary line graphs. 
The closure operation preserves the order and does not decrease the minimum degree or vertex connectivity.
Consequently, the full statement of Theorem \ref{LSZ_thm} follows from its line graph case by applying that case to $\operatorname{cl}(G)$.
As every line graph is claw-free, Theorem \ref{LSZ_thm} is therefore equivalent to its restriction to line graphs.
\medskip
	
In 1986, Thomassen \cite{T} posed a conjecture on the Hamiltonicity of line graphs, as stated in Conjecture \ref{conj1} (i). Earlier, Matthews and Sumner~\cite{MS} had proposed an apparently stronger conjecture for claw-free graphs; see Conjecture~\ref{conj1} (ii). Ryj\'{a}\v{c}ek~\cite{R} later showed that these two conjectures are in fact equivalent. Ku\v{z}el and Xiong~\cite{KX} strengthened Conjecture~\ref{conj1} (i), and Ryj\'{a}\v{c}ek and Vr\'{a}na~\cite{RP} strengthened Conjecture~\ref{conj1} (ii), both to Hamilton-connectedness. Finally, Ryj\'{a}\v{c}ek and Vr\'{a}na~\cite{RP} proved that all four statements in Conjecture~\ref{conj1} are mutually equivalent.

\begin{conjecture}\label{conj1}
{\rm(i) (Thomassen \cite{T})} Every $4$-connected line graph is Hamiltonian.
\newline\indent{\rm(ii) (Matthews and Sumner \cite{MS})} Every $4$-connected claw-free graph is Hamiltonian.
\newline\indent{\rm(iii) (Ku\v{z}el and Xiong \cite{KX})} Every $4$-connected line graph is Hamilton-connected.
\newline\indent{\rm(iv) (Ryj\'{a}\v{c}ek and Vr\'{a}na \cite{RP})} Every $4$-connected claw-free graph is Hamilton-connected.	\end{conjecture}

A \emph{hypergraph} $\mathcal{H}$ is an ordered pair $(V(\mathcal{H}),E(\mathcal{H}))$, where $V(\mathcal{H})$ is the vertex set of $\mathcal{H}$ and $E(\mathcal{H})$ is a collection of not necessarily distinct nonempty subsets of $V(\mathcal{H})$, called \emph{hyperedges} of $\mathcal{H}$. We allow singleton hyperedges, which may arise under the identifications used later. The \emph{rank} of a hypergraph $\mathcal{H}$ is the maximum cardinality of a hyperedge of $\mathcal{H}$. A hypergraph of rank at most $r$ is also referred to as an \emph{$r$-hypergraph}. A $2$-hypergraph with no singleton hyperedges corresponds to a loopless multigraph.  The \emph{line graph} $L(\mathcal{H})$ has one vertex for each hyperedge of $\mathcal{H}$, with two vertices adjacent if the corresponding hyperedges intersect. Repeated hyperedges are represented by distinct vertices.

Line graphs of $3$-hypergraphs form a larger class than ordinary line graphs and may contain induced claws. In view of the role of ordinary line graphs in Theorem \ref{LSZ_thm}, it is natural to ask whether the minimum degree result in that theorem extends to this larger class. A further motivation for studying line graphs of $3$-hypergraphs comes from their close connection with Conjecture~\ref{conj1}, established by Li, Ozeki, Ryj\'{a}\v{c}ek and Vr\'{a}na \cite{LORV} in 2020.
	
\begin{theorem}[Li, Ozeki, Ryj\'{a}\v{c}ek and Vr\'{a}na~\cite{LORV}]\label{thm3}
Conjecture~\ref{conj1} is equivalent to the statement that every $4$-connected line graph of a $3$-hypergraph is Hamilton-connected.
\end{theorem}
	
A result of Gu, Lai, and Song~\cite{GLS} implies that the line graph of every $6$-edge-connected 3-hypergraph  is Hamiltonian. Kaiser and Vr\'ana~\cite{KV2} proved that every $52$-connected line graph of a $3$-hypergraph is Hamiltonian. Ozeki and Zhang~\cite{OZ} proved that every $3$-connected line graph of a $3$-hypergraph is Hamiltonian if its domination number is at most four. Here, the \emph{domination number} of a graph is the minimum cardinality of a vertex set such that every vertex outside the set has a neighbor in it.

In this paper, we extend the line graph case of Theorem \ref{LSZ_thm} to sufficiently large line graphs of $3$-hypergraphs. Under the same $3$-connectivity assumption and the same minimum degree bound $\delta(G)\geq(n+5)/10$, we give a complete characterization of the nonhamiltonian graphs in this larger class. The exceptional family is obtained from the Petersen-based line graphs appearing in Theorem \ref{LSZ_thm} by allowing certain edge deletions. We first define this family and then state our main theorem.

\begin{definition}\label{def:family}
For an integer $t\ge3$, let $P^{(t)}$ be obtained from the Petersen graph $P$ by attaching $t$ pendant edges to each vertex, and put $Q_t=L(P^{(t)})$ (see Fig.~\ref{fig:petersen-qt}).
For each $e\in E(P)$, let $x_e$ be the vertex of $Q_t$ corresponding to $e$.

The family $\cP_t$ consists of all graphs $G$ obtained
from $Q_t$ by deleting only edges whose endpoints both belong
to $\{x_e:e\in E(P)\}$, subject to the following condition.
For each $v\in V(P)$, let $F_v$ be the subgraph of $G$
induced by the three vertices corresponding to the edges
of $P$ incident with $v$.
There is no edge $e=uv\in E(P)$ such that both $F_u$ and
$F_v$ are isomorphic to $P_3$ and $x_e$ is the middle vertex
of both paths.
%
\end{definition}
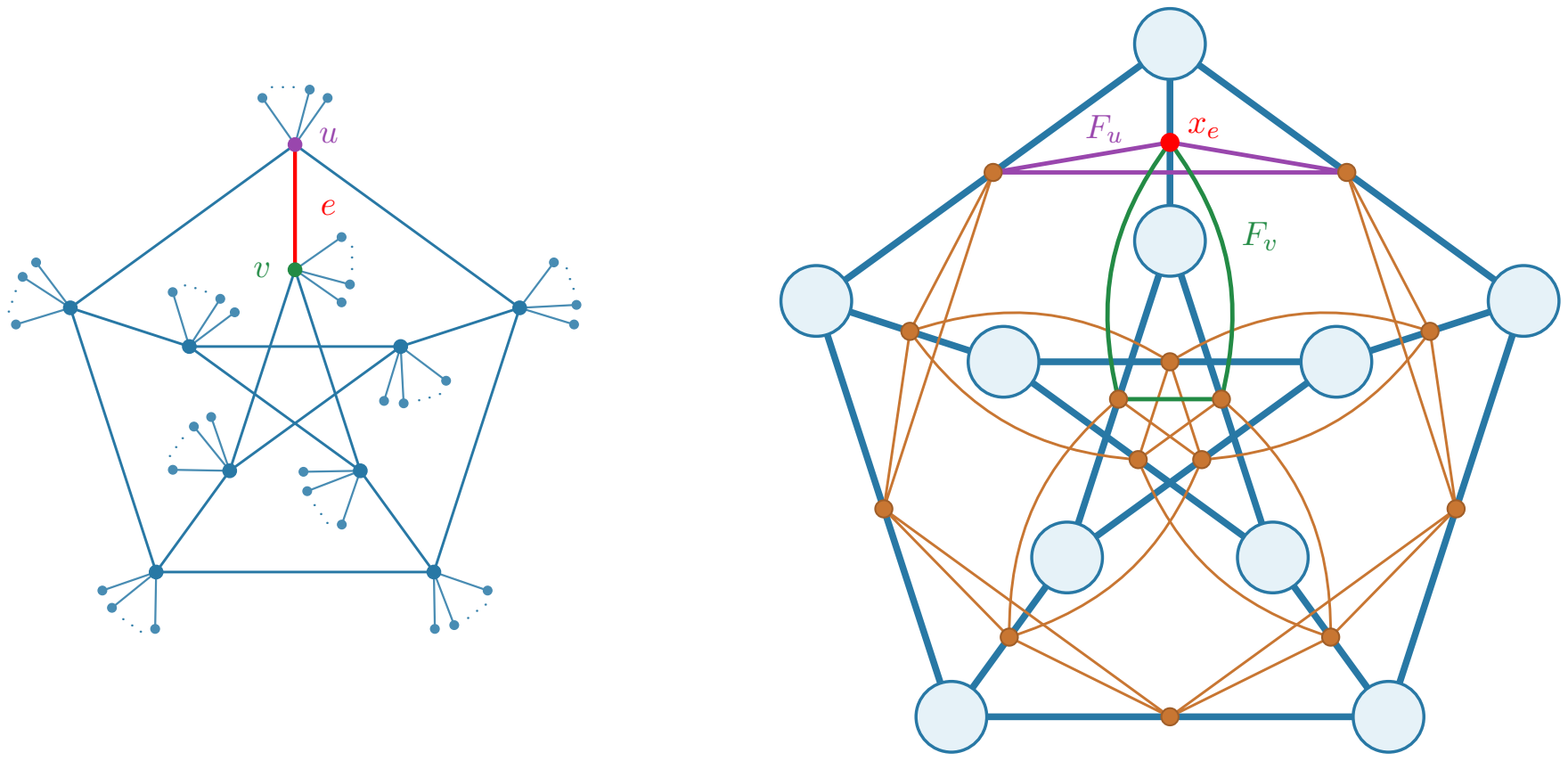
\begin{figure}[htbp]
	\centering
	\definecolor{joinblue}{HTML}{2878A5}
	\definecolor{paleblue}{HTML}{E6F2F8}
	\definecolor{edgeorange}{HTML}{C87632}
	\definecolor{Fucolor}{HTML}{9946AD}
	\definecolor{Fvcolor}{HTML}{238B45}
	\resizebox{0.70\linewidth}{!}{%
		\begin{tikzpicture}[
			fulljoin/.style={draw=joinblue,line width=2.2pt},
			localedge/.style={draw=edgeorange,line width=.85pt},
			cluster/.style={circle,draw=joinblue,fill=paleblue,
				line width=1pt,minimum size=23pt,inner sep=0pt},
			edgevertex/.style={circle,draw=edgeorange!80!black,
				fill=edgeorange,line width=.6pt,
				minimum size=5.6pt,inner sep=0pt},
			annotation/.style={font=\normalsize,text=edgeorange!85!black,
				fill=white,inner sep=1pt},
			petersenedge/.style={draw=joinblue,line width=.85pt},
			petersenvertex/.style={circle,draw=joinblue,fill=joinblue,
				minimum size=4.4pt,inner sep=0pt},
			pendantedge/.style={draw=joinblue!85,line width=.65pt},
			leafvertex/.style={circle,draw=joinblue!85,fill=joinblue!85,
				minimum size=2.8pt,inner sep=0pt},
			pendantdots/.style={font=\tiny,text=joinblue,inner sep=0pt}
			]
			
			\begin{scope}[shift={(-10,.4)}]
				\foreach \i in {0,...,4}{
					\coordinate (A\i) at ({90-72*\i}:2.7);
					\coordinate (B\i) at ({90-72*\i}:1.27);
				}
				
				\foreach \i in {0,...,4}{
					\pgfmathtruncatemacro{\j}{mod(\i+1,5)}
					\pgfmathtruncatemacro{\k}{mod(\i+2,5)}
					\draw[petersenedge] (A\i)--(A\j);
					\draw[petersenedge] (A\i)--(B\i);
					\draw[petersenedge] (B\i)--(B\k);
				}
				
				\draw[draw=red,line width=1.25pt] (A0)--(B0)
				node[midway,right=4pt,text=red] {$e$};
				
				\foreach \i in {0,...,4}{
					\pgfmathsetmacro{\outerangle}{90-72*\i}
					\pgfmathsetmacro{\innerangle}{-72*\i}
					\foreach \j/\offset in {1/-35,2/-15,3/35}{
						\coordinate (LA\i-\j)
						at ($(A\i)+({\outerangle+\offset}:.65)$);
						\coordinate (LB\i-\j)
						at ($(B\i)+({\innerangle+\offset}:.65)$);
						\draw[pendantedge] (A\i)--(LA\i-\j);
						\draw[pendantedge] (B\i)--(LB\i-\j);
						\node[leafvertex] at (LA\i-\j) {};
						\node[leafvertex] at (LB\i-\j) {};
					}
					\node[pendantdots,rotate={\outerangle-90}]
					at ($(A\i)+({\outerangle+10}:.66)$) {$\cdots$};
					\node[pendantdots,rotate={\innerangle-90}]
					at ($(B\i)+({\innerangle+10}:.66)$) {$\cdots$};
				}
				
				\foreach \i in {1,...,4}{
					\node[petersenvertex] at (A\i) {};
					\node[petersenvertex] at (B\i) {};
				}
				
				\node[petersenvertex,draw=Fucolor,fill=Fucolor]
				at (A0) {};
				\node[petersenvertex,draw=Fvcolor,fill=Fvcolor]
				at (B0) {};
				
				\node[right=4pt,yshift=3pt,text=Fucolor] at (A0) {$u$};
				\node[left=4pt,text=Fvcolor] at (B0) {$v$};
			\end{scope}
			
			\foreach \i in {0,...,4}{
				\coordinate (a\i) at ({90-72*\i}:4.25);
				\coordinate (b\i) at ({90-72*\i}:2);
			}
			
			\foreach \i in {0,...,4}{
				\pgfmathtruncatemacro{\j}{mod(\i+1,5)}
				\pgfmathtruncatemacro{\k}{mod(\i+2,5)}
				\coordinate (p\i) at ($(a\i)!.5!(a\j)$);
				\coordinate (q\i) at ($(a\i)!.5!(b\i)$);
				\coordinate (r\i) at ($(b\i)!.5!(b\k)$);
			}
			
			\foreach \i in {0,...,4}{
				\pgfmathtruncatemacro{\j}{mod(\i+1,5)}
				\pgfmathtruncatemacro{\k}{mod(\i+2,5)}
				\draw[fulljoin] (a\i)--(p\i)--(a\j);
				\draw[fulljoin] (a\i)--(q\i)--(b\i);
				\draw[fulljoin] (b\i)--(r\i)--(b\k);
			}
			
			\foreach \i in {0,...,4}{
				\pgfmathtruncatemacro{\j}{mod(\i+4,5)}
				\draw[localedge] (p\i)--(p\j)--(q\i)--cycle;
			}
			
			\foreach \i in {0,...,4}{
				\pgfmathtruncatemacro{\j}{mod(\i+3,5)}
				\draw[localedge] (q\i) to[bend left=25] (r\i);
				\draw[localedge] (q\i) to[bend right=25] (r\j);
				\draw[localedge] (r\i)--(r\j);
			}
			
			\draw[draw=Fucolor,line width=1.4pt]
			(p0)--(p4)--(q0)--cycle;
			\draw[draw=Fvcolor,line width=1.4pt]
			(q0) to[bend left=25] (r0);
			\draw[draw=Fvcolor,line width=1.4pt]
			(q0) to[bend right=25] (r3);
			\draw[draw=Fvcolor,line width=1.4pt]
			(r0)--(r3);
			
			\foreach \i in {0,...,4}{
				\node[cluster] at (a\i) {};
				\node[cluster] at (b\i) {};
				\node[edgevertex] at (p\i) {};
				\node[edgevertex] at (r\i) {};
			}
			\foreach \i in {1,...,4}{
				\node[edgevertex] at (q\i) {};
			}
			
			\node[edgevertex,draw=red,fill=red] at (q0) {};
			
			\node[annotation,text=Fucolor] at (-.75,3.28) {$F_u$};
			\node[annotation,text=Fvcolor] at (1.03,2.05) {$F_v$};
			\node[annotation,anchor=west,text=red]
			at ($(q0)+(.16,.16)$) {$x_e$};
			
		\end{tikzpicture}%
	}
	\caption{The graphs $P^{(t)}$ (left) and $Q_t=L(P^{(t)})$ (right).
Each sky-blue circle in $Q_t$ represents a clique of order $t$.
}
	\label{fig:petersen-qt}
\end{figure}

Every member of $\cP_t$ is $3$-connected and is the line graph of a $3$-hypergraph.
The latter assertion follows from Berge's characterization~\cite{B}.
	
\begin{theorem}\label{thm:main}
Let $\mathcal{H}$ be a 3-hypergraph, and let $G=L(\mathcal{H})$ be a $3$-connected graph of order $n\ge4936$. If $\delta(G)\ge(n+5)/10$, then $G$ is nonhamiltonian if and only if $G\in\cP_t$ for some integer $t\ge3$.
\end{theorem}

The constant $5$ in the numerator cannot be replaced by $4$.
Indeed, attach $t$ pendant edges to nine vertices of the Petersen graph
and $t+1$ pendant edges to the remaining vertex, where $t\ge3$.
The resulting line graph is $3$-connected and nonhamiltonian,
has order $n=10t+16$ and minimum degree $(n+4)/10$,
and belongs to none of the families $\mathcal P_s$.
These examples also show that the coefficient $1/10$ cannot be decreased.

For $n\geq4936$, Theorem \ref{thm:main} also implies Theorem \ref{LSZ_thm}. Indeed, let $F$ be a nonhamiltonian claw-free graph satisfying the hypotheses of Theorem \ref{LSZ_thm}. Applying Theorem \ref{thm:main} to $\operatorname{cl}(F)$ yields $\operatorname{cl}(F)\in\mathcal P_t$ for some $t\geq 3$. Since $Q_t$ is the only member of $\mathcal P_t$ with no eligible vertex, we have $\operatorname{cl}(F)=Q_t$, as required. Indeed, if an edge $x_ex_f$ of $Q_t$ is deleted, then $e$ and $f$ share an endpoint $v$ in $P$.
Every vertex corresponding to a pendant edge at $v$ then has a connected but noncomplete neighborhood and is therefore eligible.

Our proof adapts and simplifies the arguments of~\cite{FF,LSZ}.
In particular, it avoids the use of~\cite[Lemma~1]{FF}
and does not require the auxiliary vertex sets
$X_0$, $Y_0$, $X_I$, and $Y_I$ used in those papers;
see~\cite[pp.~302--303]{FF},
\cite[Lemma~3.2]{LSZ}, and~\cite[p.~499]{LSZ}.
This simplification enables us to extend their approach
to line graphs of $3$-hypergraphs and to prove Theorem \ref{thm:main}.

We work in the incidence graph $IG(\mathcal{H})$.  In Section~2, we introduce the notion of a $W$-dominating closed trail and establish some preliminary lemmas needed for the proof.
Section~3 proves \cref{thm:main}.

\section{Incidence graphs and reductions}
	
Throughout this section, $\mathcal{H}$ is a $3$-hypergraph having at least three hyperedges. We use the incidence graph convention of~\cite{LORV}.
The \emph{incidence graph} $IG(\mathcal H)$ is the bipartite graph
with \emph{black} vertex set $V(\mathcal H)$ and \emph{white} vertex set
$W=E(\mathcal H)$, in which a black vertex $v$ is adjacent
to a white vertex $e$ if and only if $v\in e$.
Thus every white vertex has degree at most three.
For $x\in V(\mathcal H)$, the \emph{degree} $d_{\mathcal H}(x)$
is the number of hyperedges containing $x$, counted with multiplicity;
equivalently, $d_{\mathcal H}(x)=d_{IG(\mathcal H)}(x)$.


A closed trail $T$ in the incidence graph $IG(\mathcal{H})$ is called \emph{$W$-dominating} if every vertex of $W$ has a neighbor on $T$. Equivalently, $V_b(T)$ meets every hyperedge of $\mathcal{H}$, where $V_b(T)$ is the set of black vertices contained in $T$. We allow a closed trail to consist of a single vertex. When convenient, we identify a closed trail with its underlying connected even subgraph.
Since every white vertex has degree at most three in $IG(\mathcal{H})$, 
any closed trail visits each white vertex at most once.
We use the following equivalent formulation of the characterization of Hamiltonicity in line graphs of $3$-hypergraphs due to Li, Ozeki, Ryj\'a\v{c}ek and Vr\'ana~\cite{LORV}.


\begin{theorem}[Li, Ozeki, Ryj\'a\v{c}ek and Vr\'ana]\label{thm:quasitrail}
Let $\mathcal{H}$ be a $3$-hypergraph with at least three hyperedges, and let $W$ be the white vertex set of $IG(\mathcal{H})$. The line graph $L(\mathcal{H})$ is Hamiltonian if and only if $IG(\mathcal{H})$ contains a $W$-dominating closed trail.
\end{theorem}

For distinct vertices $u,v\in V(\mathcal H)$, let
$E_{\mathcal H}(u,v)$ denote the collection of hyperedges
containing both $u$ and $v$.
The cardinality $|E_{\mathcal H}(u,v)|$ is called the
\emph{codegree} of $u$ and $v$.
We obtain the contracted hypergraph $\mathcal H/\{u,v\}$ by identifying $u$ and $v$ into a single vertex $[uv]$.
Each hyperedge $e$ meeting $\{u,v\}$ is replaced by $(e\setminus\{u,v\})\cup\{[uv]\}$, while all other hyperedges remain unchanged.
Each hyperedge retains its original label, so distinct hyperedges remain distinct even if their vertex sets become identical.
Using these labels, we regard $IG(\mathcal H)$ and $IG(\mathcal H/\{u,v\})$ as having the same white vertex set.

Under the natural correspondence between the hyperedges,
	
\begin{equation}\label{eq:monotone}
 |E(\mathcal H/\{u,v\})|=|E(\mathcal H)|, \qquad L(\mathcal H)\subseteq L(\mathcal H/\{u,v\}).
\end{equation}
	
The rank remains at most three, and the connectivity and minimum degree of the line graph do not decrease.
The following lemma shows how to lift the closed trails required by Theorem~\ref{thm:quasitrail} through this identification.
	
\begin{lemma}\label{lem:lift}
Let $\mathcal H$ be a $3$-hypergraph with two vertices $u,v$ of codegree at least three.
Let $W$ denote the common white vertex set of $IG(\mathcal H)$ and $IG(\mathcal H/\{u,v\})$.
If $IG(\mathcal H/\{u,v\})$ has a $W$-dominating closed trail,
then $IG(\mathcal H)$ has a $W$-dominating closed trail.
%
\end{lemma}
	
\begin{proof}
Let $T'$ be a $W$-dominating closed trail in $IG(\mathcal H/\{u,v\})$.
If $[uv] \notin V(T')$, take $T=T'$, which is a $W$-dominating closed trail in $IG(\mathcal H)$.

Suppose now that $[uv]\in V(T')$.
We construct a closed trail $T$ in $IG(\mathcal H)$ containing
every black vertex of $T'$ other than $[uv]$, as well as
both $u$ and $v$.
Such a trail is $W$-dominating, since $T'$ is $W$-dominating.

Choose three distinct hyperedges $e_1,e_2,e_3$ containing both $u$ and $v$, indexed so that precisely $e_1,\ldots,e_k$ have their corresponding white vertices belonging to $T'$, where $0\le k\le3$.
For each $i\le k$, the vertex $e_i$ has degree two in $T'$.
As its corresponding hyperedge has size at most two after identification,
the edge $[uv]e_i$ belongs to $T'$.
Delete these $k$ edges from $T'$.
Split $[uv]$ into $u,v$ and replace each remaining incidence edge $[uv]e$,
where $e\notin\{e_1,\ldots,e_k\}$, by an original incidence edge $ue$ or $ve$. If both are available, choose either one. Retain both $u$ and $v$, even if either is isolated, and denote the resulting subgraph by $T''$. 
%
%
The subgraph $T''$ of $IG(\mathcal H)$ contains no edges
joining $\{u,v\}$ to $\{e_1,e_2,e_3\}$.
Since $T'$ is connected, every component of $T''$ contains at least one of $u,v,e_1,\ldots,e_k$.
In $T''$,
$e_1, \dots , e_k$ have odd degree, and 
all vertices but $u,v, e_1,\ldots,e_k$ have even degree.
By the handshaking lemma, up to interchanging $u$ and $v$,
there are six possibilities for $k$ and the parities of
$d_{T''}(u)$ and $d_{T''}(v)$.
We add edges to $T''$ as follows.

\begin{itemize}
\item
If $k=0$ and 
both $u$ and $v$ have odd degree in $T''$,
then add the two edges $ue_1, ve_1$.
\item
If $k=0$ and 
both $u$ and $v$ have even degree in $T''$,
then add the four edges $ue_1, ue_2, ve_1, ve_2$.
\item
If $k=1$ and 
exactly one of $u$ and $v$ has odd degree in $T''$, say $u$,
then add the five edges $ue_1, ue_2, ve_2, ue_3, ve_3$.
\item
If $k=2$ and
both $u$ and $v$ have odd degree in $T''$,
then add the four edges $ue_1, ue_2, ue_3, ve_3$.
\item
If $k=2$ and 
both $u$ and $v$ have even degree in $T''$,
then add the four edges $ue_1, ve_2, ue_3, ve_3$.
\item
Suppose that $k=3$ and 
exactly one of $u$ and $v$ has odd degree in $T''$, say $u$.
By the handshaking lemma, the component of $T''$ containing $u$ contains another vertex of odd degree. 
Since all vertices but $u, e_1,e_2,e_3$ have even degree,
we may assume that $u$ and $e_3$ lie in the same component of $T''$.
Then add the three edges $ue_1, ve_2, ve_3$.
\end{itemize}

Let $T$ be the resulting subgraph.
These choices make the degrees of all vertices even.
The added edges also ensure that $u,v$ lie in the same component and join every other component of $T''$ to it through some $e_i$. Hence $T$ is connected and even, as desired.

\end{proof}

We also need a theorem on closed trails through prescribed vertices. The nine-point theorem of Holton, McKay, Plummer and Thomassen~\cite{HMPT} was extended by Bau and Holton~\cite{BH} to twelve prescribed vertices in a $3$-connected cubic graph, with a Petersen contraction as the exception. Chen, Lai, Li, Li and Mao~\cite{CLLLM} obtained the following version for $3$-edge-connected multigraphs.
	
A multigraph $G$ \emph{contracts to} a multigraph $J$ if $V(G)$ has a partition into connected vertex sets, called \emph{parts}, such that contracting each part to a single vertex, deleting the resulting loops, and retaining all parallel edges yields a multigraph isomorphic to $J$.
Thus, in a contraction to the Petersen graph $P$, exactly fifteen edges of $G$ have endpoints in different parts.

\begin{theorem}[Chen, Lai, Li, Li and Mao~\cite{CLLLM}]\label{thm:twelve}
Let $M$ be a $3$-edge-connected multigraph and let $S \subseteq V(M)$ with $1\le|S|\le12$. Then $M$ has a closed trail containing $S$, or $M$ contracts to the Petersen graph $P$ in such a way that every part of the corresponding partition meets $S$.
\end{theorem}
	
For our application to incidence graphs, we need the following version, which allows edge cuts of size at most two provided that they do not separate prescribed vertices.
	
\begin{lemma}\label{lem:terminal}
Let $M$ be a connected multigraph, and let $S \subseteq V(M)$ with $1\le|S|\le12$. Suppose every edge cut separating two vertices of $S$ has size at least three. Then $M$ has a closed trail containing $S$, or $M$ contracts to the Petersen graph $P$ in such a way that every part of the corresponding partition meets $S$.
\end{lemma}
\begin{proof}
Repeatedly delete a bridge $e$ together with the component
of the current graph minus $e$ that is disjoint from $S$.
Once no bridges remain, deleting any $2$-edge cut produces
exactly two connected components.
For a $2$-edge cut $\{x_1u_1,x_2u_2\}$ of the current graph,
let $U$ be the component obtained by deleting these two edges
that is disjoint from $S$, where $u_1,u_2\in U$.
Delete $U$ and the two edges $x_1u_1$ and $x_2u_2$,
and add a virtual edge $x_1x_2$ if $x_1\ne x_2$.
We call this operation a \emph{$2$-edge-cut reduction}.

These operations retain $S$ and preserve the hypothesis on
cuts separating its vertices.
Each $2$-edge-cut reduction also preserves bridgelessness.
Repeat until no $2$-edge cut remains, and apply \cref{thm:twelve}.
If it has a closed trail through $S$, undo the $2$-edge-cut
reductions in reverse order, replacing each virtual edge
$x_1x_2$ used by the trail with a $u_1u_2$-path (possibly
of length zero) in the corresponding component $U$,
together with $x_1u_1$ and $x_2u_2$.
This yields a closed trail through $S$ in $M$.

Otherwise, start from the ten parts of the Petersen contraction
and undo the $2$-edge-cut reductions in reverse order.
For each reduction, add all vertices of the corresponding
component $U$ to the part containing $x_1$.
If $x_1$ and $x_2$ lie in the same part, all restored edges
remain within that part.
If they lie in different parts, the restored edge $x_2u_2$
joins these parts in place of the virtual edge $x_1x_2$.
Finally, undo the bridge deletions in reverse order, adding
each restored component to the part containing the other
endpoint of the corresponding bridge.
Thus the restored parts give the required contraction of $M$.
\end{proof}

A \emph{vertex cover} of a hypergraph $\mathcal{H}$ is a set $S\subseteq V(\mathcal{H})$ that meets every hyperedge of $\mathcal{H}$. Such vertex covers play a central role in our proof.

\section{Proof of the main theorem}
We can now give the proof of the main theorem.
\begin{proof}[Proof of \Cref{thm:main}]
Let $\mathcal H$ be a $3$-hypergraph such that $L(\mathcal H)$
is $3$-connected, has order $n\ge4936$, and satisfies
$\delta(L(\mathcal H))\ge(n+5)/10$.
We may assume that $\mathcal H$ has no isolated vertices.

For each integer $t\ge3$, every member of $\cP_t$ is a spanning subgraph of $Q_t$.
It is easy to see that $IG(P^{(t)})$ does not contain a $W_t$-dominating closed trail,
where $W_t$ is the set of white vertices in $IG(P^{(t)})$,
and hence it follows from \cref{thm:quasitrail} that $Q_t= L(P^{(t)})$ contains no Hamiltonian cycle.
Consequently, every member of $\cP_t$ is nonhamiltonian.

For the converse, suppose that $L(\mathcal H)$ is nonhamiltonian.
Repeatedly identify two vertices whose codegree is at least three. Each step decreases the number of vertices, so the process terminates at a hypergraph $\mathcal K$ satisfying 
\begin{equation}\label{eq:codegree}
		|E_{\mathcal K}(x,y)| \le2 \qquad
\text{for each $x,y\in V(\mathcal K)$ with $x\ne y$}.
\end{equation}
%
Each identification preserves the hyperedge labels,
although the associated vertex sets may change.
Using these labels, we take the white vertex sets of
$IG(\mathcal H)$ and $IG(\mathcal K)$ to be the same set $W$.
Thus $|E(\mathcal K)|=n$, and by \eqref{eq:monotone},
$L(\mathcal H)$ is a spanning subgraph of $L(\mathcal K)$.
Consequently, $L(\mathcal K)$ is $3$-connected and
$$
\delta :=\delta(L(\mathcal{K}))\ge\delta(L(\mathcal{H})) \ge\left\lceil\frac{n+5}{10}\right\rceil\ge495.
$$
By \cref{thm:quasitrail}, $IG(\mathcal{H})$ has no $W$-dominating closed trail.
Repeated applications of Lemma~\ref{lem:lift}
therefore imply that
$IG(\mathcal K)$ has no $W$-dominating closed trail.
		
\begin{claim}\label{cl:cover}
The hypergraph $\mathcal{K}$ has a vertex cover $S$ with $1\le|S|\le12$ and $d_{\mathcal{K}}(x)\ge \lceil\delta/3\rceil+1$ for every $x\in S$.
\end{claim}
		
\begin{proof}
Let $Z=\{x\in V(\mathcal{K}):d_{\mathcal{K}}(x)\ge \lceil\delta/3\rceil+1\}$. 
If $Z$ is not a vertex cover, then some hyperedge $e$ satisfies $d_{\mathcal{K}}(x)\le\lceil\delta/3\rceil$ for every $x\in e$. However, we have
\[
d_{L(\mathcal{K})}(e) \le 3(\lceil\delta/3\rceil -1) < \delta = \delta(L(\mathcal{K})),
\]
a contradiction. Therefore, $Z$ is a vertex cover of $\mathcal{K}$. Choose an inclusion-minimal vertex cover $S\subseteq Z$, and let $s=|S|$. Since $s(\lceil\delta/3\rceil+1)\le \sum_{x\in S}d_{\mathcal{K}}(x) \le3n$ and $\lceil\delta/3\rceil+1>n/30$, we have $s<90$. 
Note that for $i=1,2,3$,
we have $i-1 \le \binom{i}{2}$.
Moreover, because $S$ meets every hyperedge and every pair of vertices has codegree at most two,
$$	
\begin{aligned}
\sum_{x\in S}d_{\mathcal{K}}(x)-n=\sum_{e\in E(\mathcal{K})}\bigl(|e\cap S|-1\bigr) \le\sum_{e\in E(\mathcal{K})}\binom{|e\cap S|}{2} =\sum_{\{x,y\}\in\binom S2}|E_{\mathcal{K}}(x,y)|\le2\binom s2=s(s-1).
\end{aligned}
$$			
			
By the inclusion-minimality of $S$, each $x\in S$ has a hyperedge $e_x$ such that $e_x\cap S=\{x\}$. If a hyperedge $f$ is adjacent to $e_x$ in $L(\mathcal{K})$ and avoids $x$, then $f$ contains a vertex of $e_x\setminus\{x\}$ and a vertex of $S\setminus\{x\}$. There are at most $2(s-1)$ such pairs, each contained in at most two hyperedges by (\ref{eq:codegree}). Hence
$$
\delta\le d_{L(\mathcal{K})}(e_x) \le d_{\mathcal{K}}(x)-1+4(s-1).
$$
Summing the resulting lower bounds on $d_{\mathcal{K}}(x)$ over $x\in S$ gives $\sum_{x\in S}d_{\mathcal{K}}(x)\ge s(\delta-4s+5).$ Combining this with $\sum_{x\in S}d_{\mathcal{K}}(x)\le n+s(s-1)$, we obtain $ s(\delta-4s+5)\le n+s(s-1), $ and therefore $ \delta s\le n+5s^2-6s. $ Suppose that $s\ge13$. Then $13\le s\le89$, and $n\le10\delta-5$ yields $ (s-10)\delta\le5s^2-6s-5. $ Consequently,
$$
\delta\le\frac{5s^2-6s-5}{s-10}=5s+44+\frac{435}{s-10}.
$$
The expression on the right is convex on $[13,89]$. Hence $\delta\le\max\left\{254, 39066/79\right\}<495,$ contrary to $\delta\ge495$. Thus $s\le12$. Since $S\subseteq Z$, every vertex of $S$ has degree at least $\lceil\delta/3\rceil+1$. \end{proof}
		
Fix a vertex cover $S$ provided by \cref{cl:cover}. We next determine the structure of $L(\mathcal{K})$.
		
\begin{claim}\label{cl:petersen}
We have $\delta(L(\mathcal{H}))=\delta=(n+5)/10$ and $L(\mathcal{K})\cong Q_{t}$, where $t = \delta-2$.
\end{claim}
\begin{proof}
Since $L(\mathcal{K})$ is connected, the incidence graph $IG(\mathcal{K})$ is connected. 
We first claim that every edge cut of $IG(\mathcal{K})$ separating two vertices of $S$ has size at least 3. Otherwise, let $R$ be the set of white endpoints of an edge cut of size at most two separating $x,y\in S$. Then $|R|\le2$, and $IG(\mathcal{K})-R$ separates $x$ from $y$. Since $d_{\mathcal{K}}(x),d_{\mathcal{K}}(y)\ge4$ and $|R|\le2$, both $x$ and $y$ still have a white neighbor in $IG(\mathcal{K})-R$.  Hence deleting the vertices of $L(\mathcal{K})$ corresponding to $R$ disconnects $L(\mathcal{K})$, contradicting its $3$-connectivity.
			
We may therefore apply Lemma~\ref{lem:terminal} to $IG(\mathcal{K})$ and $S$. Suppose that its first conclusion gives a closed trail $T$ through $S$. 
Since $S$ is a vertex cover of $\mathcal{K}$, $T$ is a $W$-dominating closed trail in $IG(\mathcal{K})$, a contradiction.
Thus, the second conclusion holds: $V(IG(\mathcal{K}))$ has a partition into ten parts $B_1,\ldots,B_{10}$, each meeting $S$, whose contraction results in $P$. Thus, each $B_i$ has three edges joining it to three distinct $B_j$'s with $j\ne i$, and there are exactly fifteen edges between different parts. For each $B_i$, the white endpoints of its three incident edges are distinct, since otherwise, deleting at most two white vertices would separate $S\cap B_i$ from $S\setminus B_i$ and, by the same argument as above, yield a vertex cut of size at most $2$ in $L(\mathcal{K})$, a contradiction. Since the edges joining $B_i$ to other $B_j$'s have distinct white endpoints, each white vertex is incident with at most one edge whose endpoints lie in different $B_i$'s.
			
Call a hyperedge \emph{internal to $B_i$} if its  white vertex and all its black neighbors lie in $B_i$. Call a hyperedge \emph{incoming to $B_i$} if its white vertex lies outside $B_i$ but has a black neighbor in $B_i$. 
Let $n_i$ be the number of white vertices in $B_i$, and let $r_i$ denote the number of hyperedges incoming to $B_i$.
Thus, we have
\begin{equation}\label{eq:branchsum}
 \sum_{i=1}^{10}n_i=n,\qquad\sum_{i=1}^{10}r_i= |E(P)| = 15.
\end{equation}

For each $i=1,\ldots,10$, choose $z_i\in S\cap B_i$.
Every hyperedge containing $z_i$ that is not internal to $B_i$
gives rise to an incidence edge
joining $B_i$ to its complement.
Since there are only three such incidence edges and
$d_{\mathcal K}(z_i)\ge 4$, at least one hyperedge containing $z_i$
is internal to $B_i$.
For a hyperedge $e$ internal to $B_i$,
every hyperedge intersecting $e$ either has its white vertex in $B_i$ or is incoming to $B_i$. Thus,

\begin{equation}\label{eq:branchdegree}
\delta\le d_{L(\mathcal{K})}(e)\le n_i-1+r_i.
\end{equation}
Choose one internal hyperedge for each part. Summing the corresponding inequalities and using \eqref{eq:branchsum}, we obtain $10\delta\le n+5$. The degree hypothesis now gives $ n+5\le10\delta(L(\mathcal{H})). $ Since $\delta(L(\mathcal{H})) \le \delta$, we have $\delta(L(\mathcal{H}))=\delta = \delta(L(\mathcal{K})) = (n+5)/10$, and equality holds in every term of the sum. Consequently, $n_i-1+r_i=\delta$ for each $i$.    
		
It follows from (\ref{eq:branchdegree}) that every hyperedge internal to $B_i$
intersects every other hyperedge whose white vertex lies
in $B_i$, as well as every hyperedge incoming to $B_i$.

Since only three incidence edges join $B_i$ to its complement, at least
\[
d_{\mathcal{K}}(z_i)-3
\ge\left\lceil\frac{\delta}{3}\right\rceil-2
\ge163
\]
hyperedges containing $z_i$ are internal to $B_i$.    We now show that every hyperedge whose white vertex lies in $B_i$, as well as every hyperedge incoming to $B_i$, contains $z_i$.
Suppose that one of these hyperedges, say $f$, does not contain $z_i$. 
By the argument above, there are at least $163$ internal hyperedges containing $z_i$, and $f$ intersects all of them. 
However, by (\ref{eq:codegree}), each vertex of $f$ belongs to at most two
hyperedges containing $z_i$.
Thus $f$ can intersect at most $2|f|\le 6$ of these hyperedges,
a contradiction.

Of the three edges joining $B_i$ to its complement, exactly $3-r_i$ have their white endpoints in $B_i$. Since these endpoints
are distinct, $B_i$ has
$$
n_i-(3-r_i)=n_i+r_i-3=\delta-2=:t
$$
internal hyperedges. Each hyperedge internal to $B_i$ contains $z_i$ and no other $z_j$. By the preceding observations, each remaining hyperedge has exactly one incidence edge joining different parts, say $B_i$ and $B_j$, and contains $z_i$ and $z_j$, but no $z_k$ with $k\notin\{i,j\}$. Thus the remaining hyperedges correspond one-to-one to the edges of $P$.

Extend the above correspondence by matching the $t$ hyperedges internal to each $B_i$ with the $t$ pendant edges at the corresponding vertex of $P$. If two hyperedges share a black vertex in $B_i$, both contain $z_i$. Thus two distinct hyperedges
intersect if and only if they contain a common $z_i$, which is equivalent to their corresponding edges in $P^{(t)}$ sharing an endpoint. Hence $L(\mathcal K)\cong L(P^{(t)})=Q_t.$
\end{proof}

By \cref{cl:petersen}, $L(\mathcal K)\cong Q_t$, where
$t=\delta-2\ge3$.
Using the isomorphism constructed in its proof and
\eqref{eq:monotone}, we regard $L(\mathcal H)$ as a spanning
subgraph of $Q_t$.
We now show that $L(\mathcal{H})$ belongs to $\mathcal P_t$.

For each $e\in E(P)$, let $x_e$ denote the vertex of $Q_t=L(P^{(t)})$ corresponding to $e$.
Under the correspondence established in the proof of Claim \ref{cl:petersen}, the vertices of $Q_t$ outside $\{x_e:e\in E(P)\}$ correspond precisely to the hyperedges internal to the parts $B_1,\ldots,B_{10}$.
Each of these vertices has degree exactly $t+2=\delta$ in $Q_t$.
Since $L(\mathcal H)$ is a spanning subgraph of $Q_t$ with minimum degree $\delta$, every edge of $Q_t$ incident with such a vertex belongs to $L(\mathcal{H})$.
Thus only edges whose endpoints both belong to
$\{x_e:e\in E(P)\}$ may be absent from $L(\mathcal{H})$.

For each $v\in V(P)$, 
let $F_v$ be the subgraph of $L(\mathcal{H})$ induced by the three vertices corresponding to the edges of $P$ incident with $v$. 
%
%
%
Suppose that, for some edge $e=uv\in E(P)$, both $F_u$ and $F_v$ are isomorphic to $P_3$, with $x_e$ as the middle vertex of both paths. The two ends of each path are nonadjacent, and no end of one path is adjacent to an end of the other, since $P$ is triangle-free. Thus these four vertices, together with $x_e$, induce a $K_{1,4}$. Therefore, $L(\mathcal{H})$ is not the line graph of a $3$-hypergraph, a contradiction.

		

Hence $L(\mathcal H)\in\cP_t$, completing the proof.
\end{proof}

\section*{Declarations}

\noindent \textbf{Conflict~of~interest}:
The authors declare that they have no known competing financial interests or personal relationships that could have appeared to influence the work reported in this paper.

\noindent \textbf{Data~availability}:
Data sharing is not applicable to this paper as no datasets were generated or analyzed during the current study.

\noindent \textbf{AI~usage}:
The authors used ChatGPT (OpenAI) to assist with language editing and the refinement of the presentation of the manuscript. The authors carefully reviewed and revised all AI-assisted text and take full responsibility for the content of the manuscript.

\end{document}